\documentclass[reqno]{amsart}
\usepackage{amssymb,hyperref}

\theoremstyle{plain}
\newtheorem{thm}{Theorem}
\newtheorem{lem}{Lemma}
\newtheorem{cor}{Corollary}

\theoremstyle{definition}
\newtheorem{ex}{Example}
\newtheorem{rem}{Remark}

\renewcommand{\Re}{\mathrm{Re}}

\title
[Starlikeness of order $\alpha$ rerated to Mocanu functions]
{Starlikeness of order $\alpha$ \\
rerated to Mocanu functions}

\author{Hitoshi Shiraishi}
\address{Hitoshi Shiraishi \newline
Department of Mathematics \newline
Kinki University \newline
Higashi-Osaka, Osaka 577-8502, Japan}
\email{shiraishi@math.kindai.ac.jp}

\author{Shigeyoshi Owa}
\address{Shigeyoshi Owa \newline
Department of Mathematics \newline
Kinki University \newline
Higashi-Osaka, Osaka 577-8502, Japan}
\email{owa@math.kindai.ac.jp}


\date{}

\begin{document}

\begin{abstract}
For analytic functions $f(z)$ in the open unit disk $\mathbb{U}$ with $f(0)=f'(0)-1=0$,
P. T. Mocanu (Mathematica (Cluj), {\bf 11}(34) (1969)) have considered Mocanu functions.
The object of the present paper is to discuss some sufficient problems for $f(z)$ to be starlike of order $\alpha$.
\end{abstract}

\begin{flushleft}
This paper was published in the journal: \\
Acta Univ. Apulensis Math. Inform. Special Issue (2009), 1037--1046.
\end{flushleft}
\hrule

\

\

\maketitle

\section{Introduction}

\

Let $\mathcal{A}_n$ denote the class of functions
$$
f(z)=z+a_{n+1}z^{n+1}+a_{n+2}z^{n+2}+ \ldots
\qquad(n=1,2,3,\ldots)
$$
that are analytic in the open unit disk $\mathbb{U}=\{z \in \mathbb{C}:|z|<1\}$
and $\mathcal{A}=\mathcal{A}_1$.

We denote by $\mathcal{S}$ the subclass of $\mathcal{A}_n$ consisting of univalent functions $f(z)$ in $\mathbb{U}$.

\

Let $\mathcal{S^{*}}(\alpha)$ be defined by
$$
\mathcal{S^{*}}(\alpha)
=\left\{f(z)\in\mathcal{A}_n:
\Re\left(\frac{zf'(z)}{f(z)}\right)>\alpha,\ 
0\leqq {}^{\exists} \alpha<1 \right\}.
$$
We denote by $\mathcal{S}^{*} = \mathcal{S}^{*}(0)$.

\

The basic tool in proving our results is the following lemma due to Jack \cite{m1ref1}
(also, due to Miller and Mocanu \cite{m1ref2}).

\

\begin{lem} \label{jack} \quad
Let the function $w(z)$ defined by
$$
w(z)=a_nz^n+a_{n+1}z^{n+1}+a_{n+2}z^{n+2}+ \ldots
\qquad(n=1,2,3,\ldots)
$$
be analytic in $\mathbb{U}$ with $w(0)=0$.
If $\left|w(z)\right|$ attains its maximum value on the circle $\left|z\right|=r$ at a point $z_{0}\in\mathbb{U}$,
then there exists a real number $k \geqq n$ such that
$$
\frac{z_{0}w'(z_{0})}{w(z_{0})}=k.
$$
\end{lem}

\

\section{Main results}

\

Applying Lemma \ref{jack},
we drive the following lemma.

\

\begin{lem} \label{d1lem4} \quad
If $f(z)\in\mathcal{A}_n$ satisfies
$$
\left|(\beta-\gamma)\frac{zf'(z)}{f(z)}+\gamma\left(1+\frac{zf''(z)}{f'(z)}\right)\right|
<\frac{1}{1+\rho}|n\rho\gamma-\beta|
\qquad(z\in\mathbb{U})
$$
for some complex $\beta$, $\gamma$
and some real $\rho>0$
such that $\Re\left(\dfrac{\beta}{\gamma}\right)<n\rho$,
then
$$
\left|\frac{f(z)}{zf'(z)}-1\right|<\rho
\qquad(z\in\mathbb{U}).
$$
\end{lem}

\

\begin{proof}\quad
Let us define $w(z)$ by
\begin{align}
w(z) &= \frac{f(z)}{zf'(z)}-1
\label{d1lem4eq1}\\
&= b_{n}z^{n} + b_{n+1}z^{n+1} + \ldots
\qquad(z\in\mathbb{U}). \nonumber
\end{align}

Then, clearly,
$w(z)$ is analytic in $\mathbb{U}$ and $w(0)=0$.
Differentiating both sides in (\ref{d1lem4eq1}),
we obtiain
$$
zw'(z)=1-\frac{f(z)}{zf'(z)}\left(1+\frac{zf''(z)}{f'(z)}\right)
\qquad(z\in\mathbb{U}),
$$
and therefore,
\begin{align*}
\left|(\beta-\gamma)\frac{zf'(z)}{f(z)}+\gamma\left(1+\frac{zf''(z)}{f'(z)}\right)\right|
&= \left|(\beta-\gamma)\frac{1}{1+w(z)} +\gamma\frac{1-zw(z)}{1+w'(z)}\right|\\
&= \left|\frac{1}{1+w(z)}\right||\gamma zw'(z)-\beta|\\
&< \frac{1}{1+\rho}|n\rho\gamma-\beta|
\qquad(z\in\mathbb{U}).
\end{align*}

If there exists a point $z_{0} \in \mathbb{U}$ such that
$$
\max_{\left| z \right| \leqq \left| z_{0} \right|} \left| w(z) \right|
= \left| w(z_{0}) \right|
= \rho,
$$
then Lemma \ref{jack} gives us that
$w(z_{0})=\rho e^{i \theta}$ and $z_{0}w'(z_{0})=kw(z_{0})$ ($k \geqq n$).

Thus we have
\begin{align*}
\left|(\beta-\gamma)\frac{z_0f'(z_0)}{f(z_0)}+\gamma\left(1+\frac{z_0f''(z_0)}{f'(z_0)}\right)\right|
&= \left|\frac{1}{1+w(z_0)}\right||\gamma z_0w'(z_0)-\beta|\\
&= \left|\frac{1}{1+w(z_0)}\right||\gamma kw(z_0)-\beta|\\
&\geqq \frac{1}{1+|w(z_0)|}|\gamma n\rho-\beta|\\
&\geqq \frac{1}{1+\rho}|n\rho\gamma-\beta|.
\end{align*}

This contradicts our condition in the lemma.
Therefore,
there is no $z_0\in\mathbb{U}$ such that $|w(z_0)|=\rho$.
This means that $|w(z)|<\rho$ for all $z\in\mathbb{U}$.
\end{proof}

\

We consider a new application for Lemma \ref{d1lem4}.

\

\begin{thm} \label{d1thm4} \quad
If $f(z)\in\mathcal{A}_n$ satisfies
$$
\left|(\beta-\gamma)\frac{zf'(z)}{f(z)}+\gamma\left(1+\frac{zf''(z)}{f'(z)}\right)\right|
<\frac{1}{2}|n\gamma-\beta|
\qquad(z\in\mathbb{U})
$$
for some real $0<\alpha\leqq\dfrac{1}{2}$,
some complex $\beta$ and $\gamma$
such that $\Re\left(\dfrac{\beta}{\gamma}\right)<n$,
or
$$
\left|(\beta-\gamma)\frac{zf'(z)}{f(z)}+\gamma\left(1+\frac{zf''(z)}{f'(z)}\right)\right|
<|n\gamma(1-\alpha)-\alpha\beta|
\qquad(z\in\mathbb{U})
$$
for some real $\dfrac{1}{2}\leqq\alpha<1$,
some complex $\beta$ and $\gamma$
such that $\Re\left(\dfrac{\beta}{\gamma}\right)<n \left(\dfrac{1}{\alpha}-1\right)$,
then
$$
\left|\frac{f(z)}{zf'(z)}-\frac{1}{2\alpha}\right|<\frac{1}{2\alpha}
\qquad(z\in\mathbb{U}).
$$

That is $f(z)\in\mathcal{S}^{*}(\alpha)$.
\end{thm}

\

\begin{proof}\quad
Putting $\rho=1$ for $0<\alpha\leqq\dfrac{1}{2}$,
or $\rho=\dfrac{1}{\alpha}-1$ for $\dfrac{1}{2}\leqq\alpha<1$,
the condition of Theorem \ref{d1thm4} satisfies Lemma \ref{d1lem4}.

Letting
$$
g(z)=\frac{f(z)}{zf'(z)}-\frac{1}{2\alpha}
\qquad(z\in\mathbb{U}),
$$
we have
$$
\frac{f(z)}{zf'(z)}-1=g(z)-\left(1-\frac{1}{2\alpha}\right)
\qquad(z\in\mathbb{U}).
$$

We want to proove $|g(z)|<\dfrac{1}{2\alpha}$ in $\mathbb{U}$.

From Lemma \ref{d1lem4},
we see
\begin{align*}
|g(z)|
& = \left| 1 -\frac{1}{2\alpha} +\left( \frac{f(z)}{zf'(z)} -1 \right) \right| \\
& = \left| 1 -\frac{1}{2\alpha} +\varepsilon e^{i\theta} \right|
\qquad(z\in\mathbb{U}\ ;\ 0\leqq\varepsilon<\rho).
\end{align*}

When $0<\alpha\leqq\dfrac{1}{2}$,
we obtain
\begin{align*}
|g(z)|
& = \left| 1 -\frac{1}{2\alpha} +\varepsilon e^{i\theta} \right| \\
&< \frac{1}{2\alpha}-1+\varepsilon \\
&< \frac{1}{2\alpha}
\qquad(z\in\mathbb{U}\ ;\ 0\leqq\varepsilon<1).
\end{align*}

When $\dfrac{1}{2}\leqq\alpha<1$,
we have
\begin{align*}
|g(z)|
& = \left| 1 -\frac{1}{2\alpha} +\varepsilon e^{i\theta} \right| \\
&< 1-\frac{1}{2\alpha}+\varepsilon \\
&< \frac{1}{2\alpha}
\qquad\left(z\in\mathbb{U}\ ;\ 0\leqq\varepsilon<\frac{1}{\alpha}-1 \right).
\end{align*}

This completes the proof of the theorem.
\end{proof}

\

Making $\beta=1$ and some real $\gamma$ in Theorem \ref{d1thm4},
we have Corollary \ref{d1cor4}.

\

\begin{cor} \label{d1cor4} \quad
If $f(z)\in\mathcal{A}_n$ satisfies
$$
\left|(1+\gamma)\frac{zf'(z)}{f(z)}-\gamma\left(1+\frac{zf''(z)}{f'(z)}\right)\right|
<\frac{1}{2}(n\gamma+1)
\qquad(z\in\mathbb{U})
$$
for some real $0<\alpha\leqq\dfrac{1}{2}$ and $\gamma$
with $\gamma>-\dfrac{1}{n}$,
or
$$
\left|(1+\gamma)\frac{zf'(z)}{f(z)}-\gamma\left(1+\frac{zf''(z)}{f'(z)}\right)\right|
<n\gamma(1-\alpha)+\alpha
\qquad(z\in\mathbb{U})
$$
for some real $\dfrac{1}{2}\leqq\alpha<1$ and $\gamma$
with $\gamma>-\dfrac{\alpha}{n(1-\alpha)}$,
then
$$
\left|\frac{f(z)}{zf'(z)}-\frac{1}{2\alpha}\right|<\frac{1}{2\alpha}
\qquad(z\in\mathbb{U}).
$$

That is $f(z)\in\mathcal{S}^{*}(\alpha)$.
\end{cor}

\

\begin{ex} \label{d1ex4} \quad
We take
\begin{align*}
f(z)
&= \left(\frac{\beta}{\gamma}\int_{0}^{z}t^{\frac{\beta}{\gamma}-1}\left(1+\frac{\bar{\beta}}{S}t^n\right)^{\frac{S^2-|\beta|^2}{n\bar{\beta}\gamma}}dt \right)^{\frac{\gamma}{\beta}}\\
&= z+a_{n+1}z^{n+1}+a_{n+2}z^{n+2}+ \ldots,
\qquad(z\in\mathbb{U})
\end{align*}
where $S=\dfrac{1}{2}|n\gamma-\beta|$
for some real $0<\alpha\leqq\dfrac{1}{2}$,
some complex $\beta$ and $\gamma$
such that $\Re\left(\dfrac{\beta}{\gamma}\right)<n$,
or $S=|n\gamma(1-\alpha)-\alpha\beta|$
for some real $\dfrac{1}{2}\leqq\alpha<1$,
some complex $\beta$ and $\gamma$
such that $\Re\left(\dfrac{\beta}{\gamma}\right)<n\left(\dfrac{1}{\alpha}-1\right)$.

Then $f(z)$ satisfies Theorem \ref{d1thm4}.
Because,
we can get
$$
\left| (\beta-\gamma)\frac{zf'(z)}{f(z)} +\gamma \left( 1+ \frac{zf''(z)}{f'(z)} \right) \right|
= \left| \frac{\gamma+Sz^n}{1+\frac{\bar{\gamma}}{S}z^n} \right|
< S
\qquad(z\in\mathbb{U}).
$$
\end{ex}

\

Also, using Lemma \ref{jack},
we get

\

\begin{rem} \label{d1re4} \quad
If the function $f(z)\in\mathcal{A}_n$ with $\dfrac{f(z)f'(z)}{z}\neq0$ in $\mathbb{U}$ satisfies
$$
\Re \left\{(1-\alpha) \frac{zf'(z)}{f(z)} + \alpha \left( 1+ \frac{zf''(z)}{f'(z)} \right) \right\}
>0
\qquad(z\in\mathbb{U})
$$
for some real constant $\alpha$,
then $f(z)$ is said to be $\alpha$-convex in $\mathbb{U}$.
We denote this class by $\mathcal{M}_\alpha$.
The class of $\alpha$-convex functions (Mocanu functions) in $\mathbb{U}$ was introduced by Mocanu \cite{d1ref3}.
\end{rem}

\

\begin{lem} \label{d1lem5} \quad
If $f(z)\in\mathcal{A}_n$ satisfies
$$
\left|\beta\left(\frac{zf'(z)}{f(z)}-1\right)+\gamma\frac{zf''(z)}{f'(z)}\right|
<\frac{\rho}{1+\rho}|\beta+\gamma(n+1)|
\qquad(z\in\mathbb{U})
$$
for some real $\rho>0$,
some complex $\beta$ and $\gamma$
such that $\Re\left(\dfrac{\beta}{\gamma}\right)>-(n+1)$,
then
$$
\left|\frac{f(z)}{zf'(z)}-1\right|<\rho
\qquad(z\in\mathbb{U}).
$$
\end{lem}

\

\begin{proof}\quad
Let us define  the function $w(z)$ by
\begin{align}
w(z) &= \frac{f(z)}{zf'(z)}-1
\label{d1lem5eq1}\\
&= b_{n}z^{n} + b_{n+1}z^{n+1} + \ldots
\qquad(z\in\mathbb{U}). \nonumber
\end{align}

Then,
we have that $w(z)$ is analytic in $\mathbb{U}$ and $w(0)=0$.
Differentiating both sides in (\ref{d1lem5eq1}) and simplifying,
we obtain
$$
zw'(z)=1-\frac{f(z)}{zf'(z)}\left(1+\frac{zf''(z)}{f'(z)}\right)
\qquad(z\in\mathbb{U}),
$$
and,
hence
\begin{align*}
\left|\beta\left(\frac{zf'(z)}{f(z)}-1\right)+\gamma\frac{zf''(z)}{f'(z)}\right|
&= \left|\beta\left(\frac{1}{1+w(z)}-1\right)+\gamma\frac{zw'(z)+w(z)}{1+w(z)}\right|\\
&= \left|\frac{w(z)}{1+w(z)}\right|\left|\beta+\gamma\left(\frac{zw'(z)}{w(z)}+1\right)\right|\\
&< \frac{\rho}{1+\rho}|\beta+\gamma(n+1)|
\qquad( z\in\mathbb{U}).
\end{align*}

If there exists a point $z_{0} \in \mathbb{U}$ such that
$$
\max_{\left| z \right| \leqq \left| z_{0} \right|} \left| w(z) \right|
= \left| w(z_{0}) \right|
= \rho,
$$
then Lemma \ref{jack} gives us that
$w(z_{0})=\rho e^{i \theta}$ and $z_{0}w'(z_{0})=kw(z_{0})$ ($k \geqq n$).

Thus we have
\begin{align*}
\left|\beta\left(\frac{z_0f'(z_0)}{f(z_0)}-1\right)+\gamma\frac{z_0f''(z_0)}{f'(z_0)}\right|
&= \left|\frac{w(z_0)}{1+w(z_0)}\right|\left|\beta+\gamma\left(\frac{z_0w'(z_0)}{w(z_0)}+1\right)\right|\\
&= \left|\frac{w(z_0)}{1+w(z_0)}\right||\beta+\gamma(k+1)|\\
&\geqq \frac{\rho}{1+\rho}|\beta+\gamma(k+1)|\\
&\geqq \frac{\rho}{1+\rho}|\beta+\gamma(n+1)|,
\end{align*}
which contradicts the condition in the lemma.
Therefore,
there is no $z_0\in\mathbb{U}$ such that $|w(z_0)|=\rho$.
This means that $|w(z)|<\rho$ for all $z\in\mathbb{U}$,
that is,
that
$$
\left| \frac{f(z)}{zf'(z)} -1 \right|
< \rho
\qquad( z\in\mathbb{U}).
$$
\end{proof}

\

Further,
from \ref{d1lem5},
we have

\

\begin{thm} \label{d1thm5} \quad
If $f(z)\in\mathcal{A}_n$ satisfies
$$
\left|\beta\left(\frac{zf'(z)}{f(z)}-1\right)+\gamma\frac{zf''(z)}{f'(z)}\right|
<\frac{1}{2}|\beta+\gamma(n+1)|
\qquad(z\in\mathbb{U})
$$
for some real $0<\alpha\leqq\dfrac{1}{2}$,
some complex $\beta$ and $\gamma$
such that $\Re\left(\dfrac{\beta}{\gamma}\right)>-(n+1)$,
or
$$
\left|\beta\left(\frac{zf'(z)}{f(z)}-1\right)+\gamma\frac{zf''(z)}{f'(z)}\right|
<(1-\alpha)|\beta+\gamma(n+1)|
\qquad(z\in\mathbb{U})
$$
for some real $\dfrac{1}{2}\leqq\alpha<1$,
some complex $\beta$ and $\gamma$
such that $\Re\left(\dfrac{\beta}{\gamma}\right)>-(n+1)$,
then
$$
\left|\frac{f(z)}{zf'(z)}-\frac{1}{2\alpha}\right|<\frac{1}{2\alpha}
\qquad(z\in\mathbb{U}),
$$
which shows that $f(z)\in\mathcal{S}^{*}(\alpha).$
\end{thm}

\

\begin{proof}\quad
We take $\rho=1$ for $0<\alpha\leqq\dfrac{1}{2}$,
or $\rho=\dfrac{1}{\alpha}-1$ for $\dfrac{1}{2}\leqq\alpha<1$,
this condition satisfies Lemma \ref{d1lem5}.

We get
\begin{align*}
|g(z)|
& = \left| 1 -\frac{1}{2\alpha} +\left( \frac{f(z)}{zf'(z)} -1 \right) \right| \\
& = \left| 1 -\frac{1}{2\alpha} +\varepsilon e^{i\theta} \right|
\qquad(z\in\mathbb{U}\ ;\ 0\leqq\varepsilon<\rho)
\end{align*}
from Lemma \ref{d1lem5}.

When $0<\alpha\leqq\dfrac{1}{2}$,
we see
\begin{align*}
|g(z)|
& = \left| 1 -\frac{1}{2\alpha} +\varepsilon e^{i\theta} \right| \\
&< \frac{1}{2\alpha}-1+\varepsilon \\
&< \frac{1}{2\alpha}
\qquad(z\in\mathbb{U}\ ;\ 0\leqq\varepsilon<1).
\end{align*}

When $\dfrac{1}{2}\leqq\alpha<1$,
we obtain
\begin{align*}
|g(z)|
& = \left| 1 -\frac{1}{2\alpha} +\varepsilon e^{i\theta} \right| \\
&< 1-\frac{1}{2\alpha}+\varepsilon \\
&< \frac{1}{2\alpha}
\qquad\left(z\in\mathbb{U}\ ;\ 0\leqq\varepsilon<\frac{1}{\alpha}-1 \right).
\end{align*}
\end{proof}

\

\begin{ex} \label{d1ex5} \quad
We consider the function $f(z)$ given by
\begin{align*}
f(z)
&= \left(\frac{\beta+\gamma}{\gamma}\int_0^z t^{\frac{\beta}{\gamma}}\exp\left(\frac{S}{n\gamma}t^n\right)dt\right)^{\frac{\gamma}{\beta+\gamma}}\\
&= z+a_{n+1}z^{n+1}+a_{n+2}z^{n+2}+ \ldots
\qquad(z\in\mathbb{U})
\end{align*}
where $S=\dfrac{1}{2}|\beta+\gamma(n+1)|$
for some real $0<\alpha\leqq\dfrac{1}{2}$,
some complex $\beta$ and $\gamma$
such that $\Re\left(\dfrac{\beta}{\gamma}\right)>-(n+1)$,
or $S=(1-\alpha)|\beta+\gamma(n+1)|$
for some real $\dfrac{1}{2}\leqq\alpha<1$,
some complex $\beta$ and $\gamma$
such that $\Re\left(\dfrac{\beta}{\gamma}\right)>-(n+1)$.

We can change $f(z)$ that
$$
\beta \left( \frac{zf'(z)}{f(z)}-1 \right) +\gamma \frac{zf''(z)}{f'(z)}=Sz^n
\qquad(z\in\mathbb{U}).
$$

So,
the function $f(z)$ satisfies Theorem \ref{d1thm5}.
\end{ex}

\


\begin{thebibliography}{}

\bibitem{m1ref1}
I. S. Jack,
{\it Functions starlike and convex of order $\alpha$},
J. London Math. Soc. {\bf 3}(1971), 469--474.

\bibitem{m1ref2}
S. S. Miller and P. T. Mocanu,
{\it Second-order differential inequalities in the complex plane},
J. Math. Anal. Appl. {\bf 65}(1978), 289--305.

\bibitem{d1ref2}
S. S. Miller and P. T. Mocanu,
{\it Univalent solutions of Briot-Bouquet differential equations},
J. Diff. Eqns. {\bf 56}(1985), 297--309.

\bibitem{d1ref3}
P. T. Mocanu,
{\it Une propri\'{e}t\'{e} de convexit\'e g\'{e}n\'{e}ralise\'{e} dans lath\'{e}orie de la repr\'{e}sentation conforme},
Mathematica (Cluj), {\bf 11}(34) (1969), 127--133.

\end{thebibliography}
\end{document}